\documentclass[11pt, twoside, a4paper]{article}

\usepackage{amsmath}     				
\usepackage{amssymb}   					
\usepackage{amsthm}     				
\usepackage{graphicx}    				

\usepackage{textcomp}    				
\usepackage[T1]{fontenc} 				
\usepackage{marvosym}    				
\usepackage[sc]{mathpazo}  	 			

\usepackage{authblk} 					
\usepackage[usenames]{xcolor}  			
\usepackage{lastpage} 					

\usepackage{enumerate}	 				
\usepackage{mathtools}					
\usepackage{bm}						
\usepackage[noadjust]{cite}			

\usepackage{hyperref} 					

\definecolor{ForestGreen}{rgb}{0.15,0.416,0.18}
\definecolor{EgyptBlue}{rgb}{0.063,0.2,0.65}
\hypersetup{
	colorlinks=true,
	linkcolor=EgyptBlue,         
   citecolor=EgyptBlue,          
   urlcolor=ForestGreen          
}

\newtheorem{theorem}{Theorem}[section]
\newtheorem{corollary}[theorem]{Corollary}
\newtheorem{lemma}[theorem]{Lemma}

\theoremstyle{definition}
\newtheorem{definition}[theorem]{Definition}
\theoremstyle{definition}
\newtheorem{remark}[theorem]{Remark}
\theoremstyle{definition}

\numberwithin{equation}{section}
\numberwithin{table}{section}
\numberwithin{figure}{section}

\title{Isospectral Dirac operators}

\author{\textbf{Yuri Ashrafyan}}
\author{\textbf{Tigran Harutyunyan}\footnote{Corresponding author. Email: hartigr@yahoo.co.uk}}

\affil{Yerevan State University, Alex Manoogian 1, Yerevan, 0025, Armenia}

\newcommand\shorttitle{Isospectral Dirac operators}

\newcommand\authorsshort{Yu. A. Ashrafyan and T. N. Harutyunyan}

\newcommand{\ejqtdelogo}{tink.png}

\makeatletter
\renewcommand{\maketitle}{\bgroup\setlength{\parindent}{0pt}

\begin{picture}(20,20)
     \put(-50,-5){\includegraphics[width=2.9truecm]{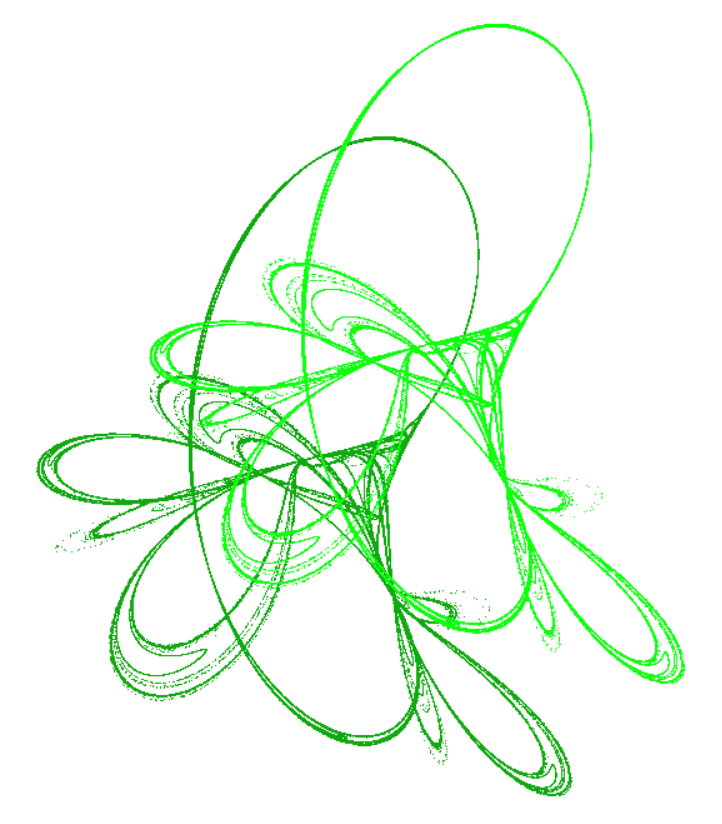}}
  \end{picture}

\vspace{1truecm}
\begin{center}{\vbox{\titlefont\@title}}\end{center}
\vspace{0.5truecm}
\begin{center}{\@author} \end{center}

\egroup
}

\renewcommand{\@fnsymbol}[1]{%
    \ifcase#1 \or {\,\Letter\!} \or $1$ \or \textasteriskcentered\or \textasteriskcentered\textasteriskcentered 
    \else\@ctrerr\fi}
\makeatother

\makeatletter
\newcommand{\hbibitem}[4]{\bibitem{#1}{#2}
\def\@tempa{#3}%
\def\@tempb{#4}%
\ifx\@tempa\@empty\ifx\@tempb\@empty{}{}\else{}{\href{https://doi.org/#4}{url}}\fi\else
{\href{http://www.ams.org/mathscinet-getitem?mr=#3}{MR#3}}\ifx\@tempb\@empty{}\else{; \href{http://dx.doi.org/#4}{url}}\fi\fi}
\makeatother

\newcommand*{\titlefont}{\fontsize{18}{21.6}\selectfont\textbf}

\makeatletter
\renewcommand\@author{\ifx\AB@affillist\AB@empty\AB@author\else
      \ifnum\value{affil}>\value{Maxaffil}\def\rlap##1{##1}%
    \AB@authlist\\affilsep]\vbox{\AB@affillist}
    \else  \AB@authors\fi\fi}
\makeatother

\makeatletter
\def\Year{2017}
\def\IssueNumber{4}
\def\ps@plain{
\def\@oddhead{\ifnum\thepage=1\hss\baselineskip8pt
\vtop to 0 pt{\vskip-0.4truecm\noindent\hbox{\hspace{1.3truecm}\Large
Electronic Journal of Qualitative Theory of Differential Equations\hss\linebreak}%
\vskip 0.1truecm
\noindent\hbox{\hspace{1.3truecm}\small \Year, No.\ {\bf \IssueNumber}, {1--\begin{NoHyper}\pageref{LastPage}\end{NoHyper}};\enspace doi:\,10.14232/ejqtde.\Year.1.\IssueNumber} \hfill \hbox{\small \href{http://www.math.u-szeged.hu/ejqtde/}{http:/\!/www.math.u-szeged.hu/ejqtde/}}
\vss}
\else
\hss\textit{\shorttitle}\hss\hbox to 0pt{\hss\thepage}\fi}\def\@oddfoot{}
\def\@evenhead{\hbox to 0pt{\thepage\hss} \hss\textit{\authorsshort}\hss}
\def\@evenfoot{}
}
\makeatother

\begin{document}

\maketitle

\pagestyle{plain}

\begin{center}
\noindent
\begin{minipage}{0.85\textwidth}\parindent=15.5pt

\medskip
\begin{center}
\noindent Received 30 December 2015, appeared 23 January 2017 
\smallskip

\noindent Communicated by Mikl\'os Horv\'ath
\end{center}
\bigskip

{\small{
\noindent {\bf Abstract.} We give the description of self-adjoint regular Dirac operators, on $[0, \pi]$, with the same spectra.}
\smallskip

\smallskip

}

\end{minipage}
\end{center}

\section{Introduction and statement of result}

Let $p$ and $q$ are real-valued, summable on $[0,\pi]$ functions, i.e.\ $p, q \in L^1_{\mathbb{R}}[0, \pi]$.
By $L(p, q, \alpha) = L(\Omega, \alpha)$ we denote the boundary-value problem for canonical Dirac system (see \cite{GaDz,GaLe,HaAz,LeSa,Ma1}):
\begin{align}
 \ell y\equiv \left\{ B \frac{d}{dx} + \Omega(x) \right\} y& =\lambda y,
\quad x \in (0, \pi),
\quad y= \begin{pmatrix}
						y_1 \\
						y_2 \\
				\end{pmatrix}
,\quad \lambda\in \mathbb{C},\label{eq1.1}\\
y_1(0)\cos\alpha+y_2(0)\sin\alpha&=0,\quad \alpha\in \left(-\frac{\pi}{2}, \frac{\pi}{2} \right], \qquad \ \qquad \qquad \qquad \label{eq1.2}\\
y_1(\pi)&=0,\label{eq1.3}
\end{align}
where
$$B=
\begin{pmatrix}
0 & 1 \\
-1 & 0 \\
\end{pmatrix}, \qquad
\Omega(x)=
\begin{pmatrix}
p(x) & q(x) \\
q(x) & -p(x) \\
\end{pmatrix}.$$
By the same $L(p, q, \alpha)$ we also denote a self-adjoint operator generated by differential expression $\ell$ in Hilbert space of two component vector-function $L^2([0, \pi]; {\mathbb{C}}^2)$ on the domain
\begin{align*}
D = {}& \bigg\{ y=
\begin{pmatrix}
y_1 \\
y_2 \\
\end{pmatrix}
;\, y_k \in AC [ 0, \pi ],\, (\ell y)_k \in L^2 [ 0, \pi ],\, k = 1, 2;\\
&\phantom{\bigg\{} y_1(0) \cos \alpha + y_2(0) \sin \alpha = 0, \quad y_1(\pi) = 0 \bigg\}
\end{align*}
where $AC [ 0, \pi ]$ is the set of absolutely continuous functions on $[0, \pi]$ (see, e.g.\ \cite{LeSa,Na}).
It is well known (see \cite{AlHrMy,GaDz,HaAz}) that under these conditions the spectra of the operator $L(p, q, \alpha)$ is purely discrete and consists of simple, real eigenvalues, which we denote by $\lambda_n= \lambda_n (p, q, \alpha)= \lambda_n (\Omega, \alpha)$, $n \in \mathbb{Z}$, to emphasize the dependence of $\lambda_n$ on quantities $p, q$ and $\alpha$.
It is also well known (see, e.g.\ \cite{AlHrMy,GaDz,HaAz}) that the eigenvalues form a sequence, unbounded below as well as above.
So we will enumerate it as $\lambda_k < \lambda_{k+1}, k \in \mathbb{Z}$, $\lambda_k>0$, when $k>0$ and $\lambda_k<0$, when $k<0$, and the nearest to zero eigenvalue we will denote by $\lambda_0$.
If there are two nearest to zero eigenvalue, then by $\lambda_0$ we will denote the negative one.
With this enumeration it is proved (see \cite{AlHrMy,GaDz,HaAz}), that the eigenvalues have the asymptotics:

\begin{equation}{\label{eq1.4}}
\lambda_n(\Omega, \alpha)= n - \frac{\alpha}{\pi} + r_n, \quad r_n = o(1), \quad n \rightarrow \pm \infty.
\end{equation}

In what follows, writing $\Omega \in A$ will mean $p, q \in A$.
If $\Omega \in L^2_{\mathbb{R}}[0, \pi]$, then we know, (see, e.g.\ \cite{HaAz}), that instead of $ r_n = o(1) $ we have:

\begin{equation}{\label{eq1.5}}
\displaystyle \sum_{n=-\infty}^{\infty} r_n^2 < \infty.
\end{equation}

Let $\varphi (x, \lambda) = \varphi(x,\lambda,\alpha, \Omega)$ be the solution of the Cauchy problem
\begin{equation}\label{eq1.6}
\ell \varphi = \lambda \varphi, \qquad \varphi (0, \lambda) = 
\begin{pmatrix}
\sin\alpha \\
-\cos\alpha \\
\end{pmatrix}
.
\end{equation}
Since the differential expression $\ell$ self-adjoint, then the components $\varphi_1(x,\lambda)$ and $\varphi_2(x,\lambda)$
of the vector-function $\varphi(x,\lambda)$ we can choose real-valued for real $\lambda$.
By $a_n=a_n(\Omega,\alpha)$ we denote the squares of the $L^2$-norm of the eigenfunctions
$\varphi_n(x,\Omega)=\varphi(x,\lambda_n(\Omega, \alpha),\alpha, \Omega)$:

\begin{equation*}
a_n=\|\varphi_n\|^2= \int_0^{\pi} |\varphi_n (x, \Omega)|^2 dx,\qquad    n \in \mathbb{Z}.
\end{equation*}
The numbers $a_n$ are called norming constants.
And by $h_n(x,\Omega)$ we will denote normalized eigenfunctions (i.e.\ $\|h_n(x)\| = 1$):

\begin{equation}\label{eq1.7}
h_n(x, \Omega) = h_n(x) = \frac{\varphi_n(x, \Omega)}{\sqrt{a_n(\Omega,\alpha)}}.
\end{equation}

It is known (see \cite{GaDz,HaAz}) that in the case of $\Omega \in L^2_{\mathbb{R}}[0, \pi]$ the norming constants have an asymptotic form:

\begin{equation}\label{eq1.8}
a_n (\Omega) = \pi + c_n, \qquad \sum_{n=-\infty}^{\infty} c_n^2 < \infty.
\end{equation}

\begin{definition}
	Two Dirac operators $L(\Omega, \alpha)$ and $L(\tilde{\Omega}, \tilde{\alpha})$ are said to be isospectral,
	if $\lambda_n(\Omega, \alpha)\hspace{-2.7pt} = \lambda_n(\tilde{\Omega}, \tilde{\alpha})$, for every $n \in \mathbb{Z}$.
\end{definition}

\begin{lemma}\label{lem1.1}
	Let $\Omega, \tilde{\Omega} \in L^1_{\mathbb{R}}[0, \pi] $ and
	the operators $L(\Omega, \alpha)$ and $L(\tilde{\Omega}, \tilde{\alpha})$ are isospectral.
	Then $\tilde{\alpha} = \alpha$.
\end{lemma}

\begin{proof}
	The proof follows from the asymptotics \eqref{eq1.4}:
	
	\begin{equation*}
	\frac{\alpha}{\pi} = \lim_{n \rightarrow \infty} (n - \lambda_n(\Omega, \alpha)) =
	\lim_{n \rightarrow \infty} (n - \lambda_n(\tilde{\Omega}, \tilde{\alpha})) = \frac{\tilde{\alpha}}{\pi}.
	\end{equation*}
\end{proof}

So, instead of isospectral operators $L(\Omega, \alpha)$ and $L(\tilde{\Omega}, \tilde{\alpha})$, we can talk about ``isospectral potentials''
$\Omega$ and $\tilde{\Omega}$.

\begin{theorem}[Uniqueness theorem]\label{thm1.1}
	The map
	\begin{equation*}
	(\Omega, \alpha) \in L^2_{\mathbb{R}}[0, \pi]\times \left( -\frac{\pi}{2}, \frac{\pi}{2} \right] \longleftrightarrow
	\{ \lambda_n(\Omega, \alpha),\, a_n(\Omega, \alpha);\, n \in \mathbb{Z} \}
	\end{equation*}
	is one-to-one.
\end{theorem}

\begin{remark}\label{rmk1.1}
	It is natural to call this a Marchenko theorem, since it is an analogue of the famous theorem of V.~A.~Marchenko \cite{Ma2}, in the case for Sturm--Liouville
	problem.
	The proof of this theorem for the case $p, q \in AC[0, \pi]$ there is in the paper \cite{Wa}.
	The detailed proof for the case $p, q \in L^2_{\mathbb{R}}[0, \pi]$ there is in \cite{Ha1} (see also \cite{FrYu,GaDz,GaLe,Ha2,Ho,WeWe}).
\end{remark}

Let us fix some $\Omega \in L^2_{\mathbb{R}}[0, \pi] $ and consider the set of all canonical potentials
\smash{$\tilde{\Omega} = 
\begin{psmallmatrix}
\tilde{p} & \tilde{q} \\
\tilde{q} & -\tilde{p} \\
\end{psmallmatrix}
$},
with the same spectra as $\Omega$:

\begin{equation*}
M^2(\Omega) = \{ \tilde{\Omega} \in L^2_{\mathbb{R}}[0, \pi]:
\lambda_n(\tilde{\Omega}, \tilde{\alpha}) = \lambda_n(\Omega, \alpha),\, n \in \mathbb{Z} \}.
\end{equation*}

Our main goal is to give the description of the set $M^2(\Omega)$ as explicit as it possible.

From the uniqueness theorem the next corollary easily follows.
\begin{corollary}\label{cor1.1}
	The map
	\begin{equation*}
	\tilde{\Omega} \in M^2 (\Omega) \leftrightarrow \{ a_n(\tilde{\Omega}),\, n \in \mathbb{Z} \}
	\end{equation*}
	is one-to-one.
\end{corollary}

Since $\tilde{\Omega} \in M^2 (\Omega)$, then $a_n(\tilde{\Omega})$ have similar to \eqref{eq1.8} asymptotics.
Since $a_n(\Omega)$ and $a_n(\tilde{\Omega})$ are positive numbers, there exist real numbers $t_n = t_n(\tilde{\Omega})$, such that
$\frac{a_n(\Omega)}{a_n(\tilde{\Omega})} = e^{t_n}$. From the latter equality and from \eqref{eq1.8} follows that
\begin{equation}\label{eq1.9}
e^{t_n} = 1 + d_n, \qquad \sum_{n=-\infty}^{\infty} d_n^2 < \infty.
\end{equation}
It is easy to see, that the sequence $\{ t_n ;\, n \in \mathbb{Z} \}$ is also from $l^2$, i.e.\ $ \sum_{n=-\infty}^{\infty} t_n^2 < \infty$.
Since all $a_n(\Omega)$ are fixed, then from the corollary {\ref{cor1.1}} and the equality
$a_n(\tilde{\Omega}) =  a_n(\Omega) e^{-t_n}$ we will get the following corollary.

\begin{corollary}\label{cor1.2}
	The map
	\begin{equation*}
	\tilde{\Omega} \in M^2 (\Omega) \leftrightarrow \{ t_n(\tilde{\Omega}),\, n \in \mathbb{Z} \} \in l^2
	\end{equation*}
	is one-to-one.
\end{corollary}

Thus, each isospectral potential is uniquely determined by a sequence $\{ t_n ;\, n \in \mathbb{Z} \}$.
Note, that the problem of description of isospectral Sturm--Liouville operators was solved in \cite{DaTr,IsMcTr,IsTr,PoTr}.

For Dirac operators the description of $M^2(\Omega)$ is given in \cite{Ha2}.
This description has a ``recurrent'' form, i.e.\ at the first in \cite{Ha2} is given the description of a family of isospectral potentials $\Omega (x, t),\ t \in \mathbb{R}$, for which only one norming constant $a_m (\Omega (\cdot, t))$ different from $a_m (\Omega)$ (namely, $a_m (\Omega (\cdot, t)) = a_m (\Omega) e^{-t}$), while the others are equal, i.e.\ $a_m (\Omega (\cdot, t)) = a_m (\Omega)$, when $n \neq m$.
\pagebreak

\begin{theorem}[\cite{Ha2}]\label{thm1.2}	Let $t \in \mathbb{R}$, $\alpha \in \big( - \frac{\pi}{2}, \frac{\pi}{2} \big]$ and
	\begin{equation*}
	\Omega(x,t) = \Omega(x) + \frac{e^{t} - 1}{\theta_m(x,t,\Omega)} \{ B h_m(x,\Omega) h_m^{*}(x,\Omega) - h_m(x,\Omega) h_m^{*}(x,\Omega) B \}
	,
	\end{equation*}
	where $\theta_n (x, t, \Omega) = 1 + (e^t - 1) \int_0^x |h_n (s, \Omega)|^2 ds$, and $*$ is a sign of transponation, e.g.\ $h_m^{*} =
		\begin{psmallmatrix}
		h_{m_1} \\
		h_{m_2} \\
		\end{psmallmatrix}^{*}
		= ( h_{m_1}, h_{m_2} )
		$.	Then, for arbitrary $t \in \mathbb{R}$, $\lambda_n(\Omega,t) = \lambda_n (\Omega)$ for all $n \in \mathbb{Z}$, $a_n(\Omega,t) = a_n (\Omega)$ for all $n \in \mathbb{Z} \backslash \{ m\}$ and $a_m(\Omega,t) = a_m (\Omega) e^{-t}$.
	The normalized eigenfunctions of the problem $L (\Omega(\cdot, t), \alpha)$ are given by the formulae:
	\begin{equation*}
	h_n(x, \Omega(\cdot, t))=\begin{dcases}
	\frac{e^{-t/2}}{\theta_m(x,t,\Omega)} h_m(x, \Omega),& \mbox{if} \ n=m, \\
	h_n(x,\Omega)-\cfrac{(e^t-1)\int_0^x h^{*}_m(s,\Omega)h_n(s,\Omega)ds}{\theta_m(x,t,\Omega)} h_m(x,\Omega), & \mbox{if} \ n \neq m.
	\end{dcases}
	\end{equation*}
\end{theorem}

Theorem \ref{thm1.2} shows that it is possible to change exactly one norming constant, keeping the others.
As examples of isospectral potentials $\Omega$ and $\tilde{\Omega}$ we can present
$ \Omega(x) \equiv 0 =
\begin{psmallmatrix}
0 & 0 \\
0 & 0 \\
\end{psmallmatrix}
$
and
\begin{equation*}
\tilde{\Omega}(x) = \Omega_{m,t}(x) = \cfrac{\pi (e^t - 1)}{\pi + (e^t - 1) x}
\begin{pmatrix}
-\sin 2 m x & \cos 2 m x \\
\cos 2 m x & \sin 2 m x \\
\end{pmatrix}
,
\end{equation*}
where $t \in \mathbb{R}$ is an arbitrary real number and $m \in \mathbb{Z}$ is an arbitrary integer.

Changing successively each $a_m (\Omega)$ by $a_m (\Omega) e^{-t_m}$, we can obtain any isospectral potential, corresponding to the sequence
$\{ t_m; m \in \mathbb{Z} \} \in l^2 $.
It follows from the uniqueness Theorem \ref{thm1.1} that the sequence, in which we change the norming constants, is not important.

In \cite{Ha2} were used the following designations:
\begin{align*}
 T_{-1} &= \{ \ldots, 0, \ldots \},\\
T_{0} &= \{ \ldots, 0, \ldots, 0, t_0,  0, \ldots, 0, \ldots \}, \\
T_{1} &= \{ \ldots, 0, \ldots, 0, 0, t_0, t_1, 0, \ldots, 0, \ldots \}, \\
T_{2} &= \{ \ldots, 0, \ldots, 0, t_{-1}, t_0, t_1, 0, \ldots, 0, \ldots \}, \\
&\vdotswithin{=}\\
T_{2n} &= \{ \ldots, 0, 0, t_{-n}, \ldots, t_{-1}, t_0, t_1, \ldots, t_{n-1}, t_{n}, 0, \ldots \}, \\
T_{2n+1} &= \{ \ldots, 0, t_{-n}, t_{-n+1}, \ldots, t_{-1}, t_0, t_1, \ldots, t_{n}, t_{n+1}, 0, \ldots \}, \\
&\vdotswithin{=}
\end{align*}
Let $\Omega(x, T_{-1}) \equiv \Omega(x)$ and

\begin{equation*}
\Omega(x, T_{m}) = \Omega(x, T_{m-1}) + \bigtriangleup \Omega(x, T_{m}), \qquad m = 0, 1, 2, \ldots,
\end{equation*}
where

\begin{equation*}
\bigtriangleup \Omega(x, T_{m}) = \cfrac{e^{t_{\tilde{m}}} - 1}{\theta_m(x, t_{\tilde{m}}, \Omega(\cdot, T_{m-1}))}
[ B h_{\tilde{m}}(x, \Omega(\cdot, T_{m-1})) h_{\tilde{m}}^{*}(\cdot)
- h_{\tilde{m}}(\cdot) h_{\tilde{m}}^{*}(\cdot) B ],
\end{equation*}
where $\tilde{m} = \frac{m+1}{2}$, if $m$ is odd and $\tilde{m} = - \frac{m}{2}$, if $m$ is even.
The arguments in others $h_{\tilde{m}}(\cdot)$ and $h_{\tilde{m}}^{*}(\cdot)$ are the same as in the first.
And after that in \cite{Ha2} the following theorem was proved.

\begin{theorem}[\cite{Ha2}]\label{thm1.3}
	Let $T = \{ t_n, n \in \mathbb{Z} \} \in l^2 $ and $\Omega \in L^2_{\mathbb{R}}[0, \pi]$. Then
	
	\begin{equation}\label{eq1.10}
	\Omega(x, T) \equiv \Omega(x) + \sum_{m=0}^{\infty} \bigtriangleup \Omega(x, T_{m}) \in M^2(\Omega).
	\end{equation}
\end{theorem}

We see, that each potential matrix $\bigtriangleup \Omega(x, T_{m})$ defined by normalized eigenfunctions \linebreak  $h_{\tilde{m}}(x, \Omega(x, T_{m-1}))$ of the previous operator $L(\Omega(\cdot, T_{m-1}), \alpha)$.
This approach we call ``recurrent'' description.

In this paper, we want to give a description of the set $M^2(\Omega)$ only in terms of eigenfunctions $h_n(x, \Omega)$ of the initial operator $L(\Omega, \alpha)$ and sequence $T \in l^2$.
With this aim, let us denote by $N(T_m)$ the set of the positions of the numbers in $T_m$, which are not necessary zero, i.e.
\begin{align*}
N(T_0) &= \{0\},\\
N(T_1) &= \{ 0, 1\}, \\
N(T_2) &= \{ -1, 0, 1 \} , \\
&\vdotswithin{=}\\
N(T_{2n}) &= \{-n, -(n-1), \ldots, 0, \ldots, n-1, n\}, \\
N(T_{2n+1})& = \{-n, -(n-1), \ldots, 0, \ldots, n, n+1\}, \\
 &\vdotswithin{=}
\end{align*}
in particular $N(T) \equiv \mathbb{Z}$.
By $S(x, T_m)$ we denote the $(m+1)\times (m+1)$ square matrix
\begin{equation}\label{eq1.11}
S(x, T_m) = \left( \delta_{ij} + (e^{t_j} - 1 ) \int_{0}^{x} h_i^{*} (s) h_j  (s) ds \right)_{i, j \in N(T_m)}
\end{equation}
where $\delta_{ij}$ is a Kronecker symbol. By $S_{p}^{(k)}(x, T_m)$ we denote a matrix which is obtained from the matrix $S(\hspace{-1pt}x, T_m\hspace{-1pt})$ by replacing the $k$th column of $S(\hspace{-1pt}x, T_m\hspace{-1pt})$ by $H_p(x, T_m) \!=\! \{\hspace{-1pt}- ( e^{t_k} \hspace{-1pt}- \hspace{-1pt}1 ) h_{k_p}\hspace{-1pt}(x)\hspace{-1pt} \}_{k \in N(T_m)}$ column, $p = 1, 2$,
Now we can formulate our result as follows.

\begin{theorem}\label{thm1.4}
	Let $T = \{ t_k \}_{k \in \mathbb{Z}} \in l^2 $ and $\Omega \in L^2_{\mathbb{R}}[0, \pi]$.
	Then the isospectral potential from $M^2(\Omega)$, corresponding to $T$, is given by the formula
	
	\begin{equation}\label{eq1.12}
	\Omega(x, T) =  \Omega(x) + G(x, x, T) B - B G(x, x, T) =  
	\begin{pmatrix}
	p(x, T) & q(x, T) \\
	q(x, T) & -p(x, T) \\
	\end{pmatrix},
	\end{equation}
	where
		\begin{equation*}
	G(x, x, T) = \cfrac{1}{\det S(x, T)}
	\sum_{k \in \mathbb{Z}}
	\begin{pmatrix}
	\det S_{1}^{(k)}(x, T) \\
	\det S_{2}^{(k)}(x, T) \\
	\end{pmatrix}
	h_k^{*}(x),
	\end{equation*}
	and $\det S(x, T) =  \lim_{m \rightarrow \infty} \det S(x, T_m)$ (the same for $\det S_p^k(x, T),\ p=1, 2$).
	
	In addition, for $ p(x, T)$ and  $q(x, T)$ we get  explicit representations:
	\begin{align*}
	p(x, T)  &= p(x) - \cfrac{1}{\det S(x, T)} \displaystyle \sum_{k \in \mathbb{Z}} \displaystyle \sum_{p=1}^2 \det S_{p}^{(k)}(x, T) h_{k_{(3-p)}}(x),\\[2mm]
	q(x, T)  &= q(x) + \cfrac{1}{\det S(x, T)} \displaystyle \sum_{k \in \mathbb{Z}} \displaystyle \sum_{p=1}^2 (-1)^{p-1} S_{p}^{(k)}(x, T) h_{k_p}(x).
	\end{align*}
\end{theorem}

\section{Proof of  Theorem \ref{thm1.4}}

The spectral function of an operator $L(\Omega, \alpha)$ defined as

\begin{equation*}
\rho(\lambda) = \begin{cases} 
\displaystyle \sum_{0 < \lambda_n \leq \lambda} \cfrac{1}{a_n(\Omega)}, \quad & \lambda > 0,   \\
- \displaystyle \sum_{\lambda < \lambda_n \leq 0} \cfrac{1}{a_n(\Omega)}, \quad& \lambda < 0,
\end{cases}
\end{equation*}
i.e.\ $\rho(\lambda)$ is left-continuous, step function with jumps in points $\lambda = \lambda_n$ equals $\frac{1}{a_n}$ and $\rho(0) = 0$.

Let $\Omega, \tilde{\Omega} \in L^2_{\mathbb{R}}[0, \pi] $ and they are isospectral.
It is known (see \cite{AlHrMy,Ar3,GaLe,LeSa}), that there exists a function $G(x, y)$ such that:

\begin{equation}\label{eq2.1}
\varphi(x,\lambda, \alpha, \tilde{\Omega}) = \varphi (x,\lambda, \alpha, \Omega) +
\int_0^x G(x, s) \varphi(s, \lambda, \alpha, \Omega) dt.
\end{equation}
It is also known (see, e.g.\ \cite{AlHrMy,GaLe,LeSa}), that the function $G(x, y)$ satisfies to the Gelfand--Levitan integral equation:

\begin{equation}\label{eq2.2}
G(x,y)+F(x,y)+ \int^x_0 G(x,s)F(s,y)ds=0,\qquad 0\leq y \leq x,
\end{equation}
where

\begin{equation}\label{eq2.3}
F(x,y)= \int_{-\infty}^{\infty} \varphi(x,\lambda, \alpha, \Omega) \varphi^{*} (y,\lambda, \alpha, \Omega)
d[\tilde{\rho}(\lambda) -\rho(\lambda)].
\end{equation}

If the potential $\tilde{\Omega}$ from $M^2(\Omega)$ is such that only finite norming constants of the operator $L(\tilde{\Omega}, \alpha)$ are different from the norming constants of the operator $L(\Omega, \alpha)$, i.e.\ $a_n(\tilde{\Omega}) = a_n (\Omega) e^{-t_n},\ n \in N(T_m)$ and the others are equal, then it means, that

\begin{equation}\label{eq2.4}
d \tilde{\rho} (\lambda) - d \rho (\lambda) =
\displaystyle \sum_{k \in N(T_m)} \left(\frac{1}{\tilde{a_k}} - \frac{1}{a_k} \right) \delta(\lambda - \lambda_k) d \lambda =
\displaystyle \sum_{k \in N(T_m)} \left(\frac{e^{t_k}-1}{a_k} \right) \delta(\lambda - \lambda_k) d \lambda,
\end{equation}
where $\delta$ is Dirac $\delta$-function. In this case the kernel $F(x, y)$ can be written in a form of a finite sum (using notation \eqref{eq1.7}):

\begin{equation}\label{eq2.5}
F(x, y) = F(x, y, T_m)= \displaystyle \sum_{k \in N(T_m)} ( e^{t_k} - 1 ) h_k(x, \Omega) h_k^{*} (y, \Omega),
\end{equation}
and consequently, the integral equation \eqref{eq2.2} becomes to an integral equation with degenerated kernel, i.e.\ it becomes to a system of linear equations and we will look for the solution in the following form:

\begin{equation}\label{eq2.6}
G(x, y, T_m)= \displaystyle \sum_{k \in N(T_m)} g_k(x) h_k^{*} (y),
\end{equation}
where $g_k(x) = 
\begin{psmallmatrix}
g_{k_1} (x) \\
g_{k_2} (x) \\
\end{psmallmatrix}
$
is an unknown vector-function.
Putting the expressions \eqref{eq2.5} and \eqref{eq2.6} into the integral equation \eqref{eq2.2} we will obtain
a system of algebraic equations for determining the functions $g_k(x)$:

\begin{equation}\label{eq2.7}
g_k(x) + \displaystyle \sum_{i \in N(T_m)} s_{ik}(x) g_i(x) = - ( e^{t_k} - 1 ) h_k(x), \qquad k \in N(T_m),
\end{equation}
where

\begin{equation*}
s_{ik}(x)= (e^{t_k} - 1 ) \int_{0}^{x} h_i^{*} (s) h_k  (s) ds.
\end{equation*}
It would be better if we consider the equations \eqref{eq2.7} for the vectors
$g_k = 
\begin{psmallmatrix}
g_{k_1} \\
g_{k_2} \\
\end{psmallmatrix}
$
by coordinates $g_{k_1}$ and $g_{k_2}$ to be a system of scalar linear equations:

\begin{equation}\label{eq2.8}
g_{k_p}(x) + \displaystyle \sum_{i \in N(T_m)}  s_{ik}(x) g_{i_p}(x) = - ( e^{t_k} - 1 ) h_{k_p}(x), \qquad k \in N(T_m), \quad p= 1, 2.
\end{equation}
The systems \eqref{eq2.8} might be written in matrix form

\begin{equation}\label{eq2.9}
S(x, T_m) g_p(x, T_m) = H_p(x, T_m), \qquad p=1, 2,
\end{equation}
where the column vectors $g_p(x, T_m) = \{ g_{k_p} (x, T_m)\}_{k \in N(T_m)},\ p = 1, 2$, and the solution can be found in the form (Cramer's rule):

\begin{equation*}
g_{k_p}(x, T_m) = \cfrac{\det S_{p}^{(k)}(x, T_m)}{\det S(x, T_m)}, \qquad k \in N(T_m), \quad p = 1, 2.
\end{equation*}

Thus we have obtained for $g_k(x)$ the following representation:
\begin{equation}\label{eq2.10}
g_k(x, T_m) = \cfrac{1}{\det S(x, T_m)}
\begin{pmatrix}
\det S_{1}^{(k)}(x, T_m) \\
\det S_{2}^{(k)}(x, T_m) \\
\end{pmatrix}
\end{equation}
and then by putting \eqref{eq2.10} into \eqref{eq2.6} we find the $G(x, y, T_m)$ function.
If the potential $\Omega$ is from $L^1_{\mathbb{R}}$, then such is also the kernel $G(x, x, T_m)$ (see \cite{Ha2}), and the relation between them gives as follows:

\begin{equation}\label{eq2.11}
\Omega(x, T_m) = \Omega(x) + G(x, x, T_m) B - B G(x, x, T_m).
\end{equation}
On the other hand we have
\begin{equation}\label{eq2.12}
\Omega(x, T_m) = \Omega(x) + \displaystyle \sum_{k=0}^{m} \bigtriangleup \Omega(x, T_{k}).
\end{equation}
So, using the Theorem \ref{thm1.3} and the equality \eqref{eq2.12} we can pass to the limit in \eqref{eq2.11}, when $m \rightarrow \infty$:
\begin{equation}\label{eq2.13}
\Omega(x, T) = \Omega(x) + G(x, x, T) B - B G(x, x, T).
\end{equation}
The potentials $\Omega(x, T)$ in \eqref{eq1.10} and \eqref{eq2.13} have the same spectral data $\{ \lambda_n(T), a_n(T) \}_{n \in \mathbb{Z}}$, and therefore they are the same and $\Omega(\cdot, T)$ defined by \eqref{eq2.13} is also from $M^2(\Omega)$.

Using \eqref{eq2.6} and \eqref{eq2.10} we calculate the expression $G(x, x, T_m) B - B G(x, x, T_m)$ and pass to the limit, obtaining for the $p(x, T)$ and $q(x, T)$ the representations:

\begin{align*}
p(x, T)  &= p(x) - \cfrac{1}{\det S(x, T)} \displaystyle \sum_{k \in N(T)} \displaystyle \sum_{p=1}^2 \det S_{p}^{(k)}(x, T) h_{k_{(3-p)}}(x),\\[2mm]
q(x, T) & = q(x) + \cfrac{1}{\det S(x, T)} \displaystyle \sum_{k \in N(T)} \displaystyle \sum_{p=1}^2 (-1)^{p-1} S_{p}^{(k)}(x, T) h_{k_p}(x).
\end{align*}
Theorem \ref{thm1.4} is proved.

For example, when we change just one norming constant (e.g.\ for $T_0$) we get two independent linear equations:
\begin{align*}
( 1 + s_{00}(x) ) g_{0_1}(x)& = - ( e^{t_0} - 1 ) h_{0_1}(x),\\[2mm]
( 1 + s_{00}(x) ) g_{0_2}(x)& = - ( e^{t_0} - 1 ) h_{0_2}(x).
\end{align*}
For the solutions we get:
\begin{align*}
g_{0_1}(x) & = - \frac{( e^{t_0} - 1 ) h_{0_1}(x)}{1 + s_{00}(x)},\\[2mm]
g_{0_2}(x)  &= - \frac{( e^{t_0} - 1 ) h_{0_2}(x)}{1 + s_{00}(x)},
\end{align*}
and for the potentials $p(x, T_0)$ and $q(x, T_0)$:
\begin{align*}
p(x, T_0)  &= p(x) + \frac{ e^{t_0} - 1 }{1 + s_{00}(x)} ( 2 h_{0_1}(x) h_{0_2}(x) ),\\[2mm]
q(x, T_0) & = q(x) + \frac{ e^{t_0} - 1 }{1 + s_{00}(x)} ( h_{0_2}^2(x) - h_{0_1}^2(x) ).
\end{align*}

\section*{Acknowledgements}
This work was supported by the RA MES State Committee of Science, in the frames of the research project No.15T-1A392.




\end{document}